\documentclass[a4paper,12pt,reqno]{amsart}
\usepackage{latexsym}
\usepackage{amsmath,amsthm,amssymb,amscd}
\usepackage{fullpage}
\usepackage{xcolor}
\usepackage{parskip}
\usepackage{mathtools}
\mathtoolsset{showonlyrefs}

\usepackage[activate={true,nocompatibility},final,tracking=true,kerning=true,spacing=true]{microtype}
\usepackage[pdfencoding=unicode,pdftex,bookmarks=true,bookmarksnumbered]{hyperref}
\definecolor{citecolour}{rgb}{0.0, 0.0, 0.8}
\colorlet{linkcolour}{green!50!black}
\hypersetup{colorlinks,breaklinks,
		linkcolor=linkcolour,citecolor=citecolour,
		filecolor=linkcolour, urlcolor=linkcolour}
\usepackage{version}

\newtheorem{prevtheorem}{Theorem}

\newtheorem{theorem}{Theorem}
\newtheorem{lemma}[theorem]{Lemma}

\newtheorem{corollary}[theorem]{Corollary}

\theoremstyle{definition}

\theoremstyle{remark}

\numberwithin{theorem}{section}
\numberwithin{equation}{section}
\numberwithin{table}{section}

\usepackage[shortlabels]{enumitem}
\begingroup
    \makeatletter
    \@for\theoremstyle:=definition,remark,plain\do{%
        \expandafter\g@addto@macro\csname th@\theoremstyle\endcsname{%
            \addtolength\thm@preskip\parskip
            }%
        }

\DeclareMathOperator{\F}{\mathfrak{F}}

\DeclareMathOperator{\Zent}{\mathbf{Z}}

\DeclareMathOperator{\cent}{\mathbf{C}}
\DeclareMathOperator{\N}{\mathfrak{N}}
\DeclareMathOperator{\norm}{\mathbf{N}}

\DeclareMathOperator{\fratt}{\mathbf{\Phi}}
\DeclareMathOperator{\aut}{\mathbf{Aut}}
\DeclareMathOperator{\hol}{\mathbf{Hol}}

\newcommand{\R}{\mathbb{\R}}

\renewcommand{\leq}{\leqslant}
\renewcommand{\geq}{\geqslant}

\newenvironment{proofof}{{\bf {Proof.} }}{\hfill $\blacksquare$ \\} 

\newenvironment{proofofthmA}{{\bf {Proof of Theorem~\ref{Thm:thmA}.} }}{\hfill $\blacksquare$ \\}
\newenvironment{proofofthmB}{{\bf {Proof of Theorem~\ref{Thm:thmB}.} }}{\hfill $\blacksquare$ \\}  
\newenvironment{proofofcor}{{\bf {Proof of Corollary~\ref{Cor:bound}.} }}{\hfill $\blacksquare$ \\}

\begin{document}

\title{A generalisation of Schenkman's theorem} 
\author{Stefanos Aivazidis$^{\dagger}$, Ina N. Safonova$^*$ and Alexander N. Skiba$^\sharp$}
\address{$^{\dagger}$Stockholm, Sweden.}
\email{stefanosaivazidis@gmail.com}
\address{$^*$Department of Mechanics and Mathematics, Belarusian State University, Minsk 220030, Belarus.}
\email{safonova@bsu.by}
\address{$^\sharp$Department of Mathematics and Technologies of Programming,
 Francisk Skorina Gomel State University, Gomel 246019, Belarus.}
\email{alexander.skiba49@gmail.com}



\begin{abstract}    
Let $G$ be a finite group and let $\F$ be a hereditary saturated formation. We denote by $\Zent_{\F}(G)$ the product of all normal subgroups $N$ of $G$ such that
 every chief 
factor $H/K$ of $G$ below $N$ is \emph{$\F$-central} in $G$, that is,  
\[
(H/K)\rtimes(G/\cent_{G}(H/K)) \in \F.
\]

A subgroup $A \leq G$ is said to be \emph{$\F$-subnormal in the sense of Kegel}, or
 \emph{$K$-$\F$-subnormal} in $G$, 
if there is a subgroup chain 
\[
A = A_0 \leq A_1 \leq \ldots \leq A_n = G
\]
such that either $A_{i-1} \trianglelefteq A_{i}$ or 
$A_i / (A_{i-1})_{A_i} \in \F$ for all $i = 1, \ldots , n$.

In this paper, we prove the following generalisation of Schenkman's Theorem on
 the centraliser of the nilpotent residual of a subnormal subgroup:
\emph{Let $\F$ be a hereditary saturated formation and let $S$ be a $K$-$\F$-subnormal subgroup of $G$. 
If $\Zent_{\F}(E) = 1$ for every subgroup $E$ of $G$ such that $S \leq E$
 then $\cent_{G}(D) \leq D$, where
$D = S^{\F}$ is the $\F$-residual of $S$.}
\end{abstract}

\maketitle

\section{Introduction}\label{Sec:Intro}
Throughout this paper, all groups are assumed to be finite and $G$ will always denote
 a finite group. 

Schenkman~\cite{schenk} proved that if $S$ is a subnormal subgroup of the group $G$
 and if $\cent_G(S)=1$, 
then $\cent_G(S^{\N}) \leq S^{\N}$, where $S^{\N}$ is the nilpotent residual of $S$. 
A particular case of Schenkman's Theorem is $G=S$: 
if $G$ is a group with trivial centre then $G^{\N} \geq \cent_G(G^{\N})$. 
A normal subgroup of a group $G$ which contains its own centraliser in $G$ we call
 \emph{large}. 
Thus the particular case of Schenkman's Theorem we mentioned can be restated thus:
\emph{in a centreless group the nilpotent residual is a large subgroup.}                                                                          

Our goal in this note is to state and prove a generalisation of Schenkman's Theorem 
valid for all hereditary saturated formations.
 
Before continuing, we need some concepts from the theory of formations. 
Recall that $G^{\F}$ denotes the \emph{$\F$-residual} of $G$, that is, the 
intersection of all normal subgroups $N$ of $G$ with $G/N \in \F$. 
A non-empty class of groups $\F$ is said to be a \emph{formation} if every
homomorphic image of $G/G^{\F}$ belongs to $\F$ for every group $G$.
A formation $\F$ is called \emph{saturated} if $G\in \F$ whenever $G^{\F} \leq \fratt(G)$
and \emph{hereditary} (A.I. Mal'cev \cite{malcevbook})
 if $H \in \F$ whenever $H \leq G \in \F$.

Now let $R/S$ be a normal section of $G$ (meaning that $S, R \trianglelefteq G$) 
and let $K$ be a normal subgroup of $G$ such that $K \leq \cent_{G}(R/S)$. 
Then we can form the semidirect product 
\[
[R/S](G/K),
\]
where 
\[
(rS)^{gK}=g^{-1}rgK
\]
for all $rS\in R/S$ and $gK\in G/K$. Moreover, we will suppress the
$\rtimes$ symbol when the implied action is the one above.
The section $R/S$ is \emph{$\F$-central} in $G$~\cite{skibabook} 
if for some normal subgroup $K$ of $G$ such that  $K\leq \cent_{G}(R/S)$ we have 
$[R/S] (G/K) \in \F$; otherwise $R/S$ is called \emph{$\F$-eccentric} in $G$. 
If $R/1$ is $\F$-central in $G$ then we say that $R$ is $\F$-central in $G$. 

The importance of these concepts stems from a result of Barnes and Kegel~\cite{barnes} which asserts that if  
$\F$ is a formation then every chief factor of a group $G \in \F$ is $\F$-central in $G$. 
Moreover, if $\F$ is any non-empty saturated formation and every chief factor of $G$ is $\F$-central in $G$ 
then $G\in \F$ (see Lemma~\ref{Lem:Fifth} below).
 
A subgroup $A \leq G$ is said to be \emph{$\F$-subnormal in the sense
 of Kegel}~\cite{kegel78} 
or \emph{$K$-$\F$-subnormal} in $G$~\cite[Defn. 6.1.4]{classes} if there is a subgroup chain
\[
A=A_0 \leq A_1 \leq \ldots \leq A_n = G
\]
such that either $A_{i-1} \trianglelefteq A_i$ or $A_i/(A_{i-1})_{A_i} \in \F$ for all $i=1, \ldots , n$.

Note that the $K$-$\F$-subnormal subgroups play an important role in many branches of formation theory
and their study is related to the investigations of many authors (cf. \cite[Chap. 6]{classes}). 

The main result of the work reported here is the following.  

\begin{prevtheorem}\label{Thm:thmA}
Let $\F$ be a hereditary saturated formation. 
Let  $S$ be a $K$-$\F$-subnormal subgroup of $G$ and suppose that 
$\Zent_{\F}(E)=1$ for every subgroup $E$ such that $S\leq E$. 
Then $\cent_{G}(S^{{\F}})\leq S^{{\F}}$.
\end{prevtheorem}

In this theorem, $\Zent_{\F}(E)$ denotes the \emph{$\F$-hypercentre} of $E$
which is defined as the product of all normal subgroups $N$ of $E$ such that either
$N=1$ or $N > 1$ and every chief factor of $E$ below $N$ is $\F$-central in $E$. 

Observe that if for a subgroup $S$ of $G$ we have $\cent_G(S) = 1$ then
$\Zent_{\N}(E)=\Zent_{\infty}(E)=1$ for every subgroup $E$ of $G$ such 
that $S\leq E$, where $\mathfrak{N}$ is the class of all nilpotent groups.
Note, also, that every subnormal subgroup is $K$-$\F$-subnormal. 
Therefore, Schenkman's original theorem is a direct consequence of Theorem~\ref{Thm:thmA}. 

\begin{corollary}[{cf. \cite{schenk} or \cite[Thm. 9.21]{isaacs}}]
Let $S$ be a subnormal subgroup of $G$ and suppose that $\cent_G(S) = 1$. 
Then $\cent_{G}(S^{\N}) \leq S^{\N}$.
\end{corollary}

We discuss some further applications of Theorem~\ref{Thm:thmA} in Section~\ref{Sec:Final}. 

The following result, which we record here separately, is an important step in the proof of Theorem~\ref{Thm:thmA}. 

\begin{prevtheorem}\label{Thm:thmB} 
Let $\F$ be a saturated formation. 
If $G$ has no non-trivial normal $\F$-central subgroups then 
$G^{\F}$ is a large subgroup of $G$, i.e. $\cent_{G}(G^{\F})\leq G^{\F}$.
\end{prevtheorem}

\begin{corollary}\label{Cor:bound}
Let $\F$ be a saturated formation and 
let $U$ be the $\F$-residual of $G$. 
If $U \cap \Zent_{\F}(G) = 1$ then 
\begin{equation}
\left\lvert G/\Zent_{\F}(G) \right\rvert \leq \left\lvert \hol(U) \right\rvert, 
\end{equation}
where $\hol(U) = U \rtimes \aut(U)$ is the familiar holomorph of the group $U$.
\end{corollary}


We prove our results in the next section and our proofs are non-standard. In 
fact, we are developing here a novel method for proving results within the context of 
formation theory that does not use the complex machinery of the theory. 
This should make the present article accessible to a wide audience.

\section{Auxiliary results and Proof of Theorems~\ref{Thm:thmA} and~\ref{Thm:thmB}}\label{Sec:Proof}

Let $D = M \rtimes A$ and $R = N \rtimes B$. Then the pairs  $(M,A)$ and $(R,B)$ are said to be 
\textbf{equivalent} provided there are isomorphisms $f : M \to N$ and 
$g : A \to B$ such that 
\[
f(a^{-1}ma)=g(a^{-1})f(m)g(a)
\] 
for all $m \in M$ and $a \in A$. 

In fact, the following lemma is known (cf. \cite[Lemma 3.27]{skibabook}).
 
\begin{lemma}\label{Lem:First}
Let  $D = M \rtimes A$ and $R = N \rtimes B$. 
If the pairs $(M, A)$ and $(R,B)$ are equivalent then $D \cong R$.
\end{lemma}

\begin{lemma}\label{Lem:Second}
Let $\F$ be a hereditary formation and let $R/S$ be an $\F$-central normal section of $G$ and 
$K \leq L$  normal subgroups of $G$ such that $L \leq \cent_{G}(R/S)$ and $[R/S](G/K) \in \F$.  
Then the following statements hold.
\begin{enumerate}[label={\upshape(\roman*)}]
\item\label{Lem:Seconda} $[R/S](G/L) \in \F$.
\item\label{Lem:Secondb} If $E$ is a subgroup of $G$ then $(E\cap R)\big/(E\cap S)$ is $\F$-central in $E$.
\item\label{Lem:Secondc} $[T/S] (G/K) \in \F$ and  $[R/T](G/K) \in \F$ 
for every normal subgroup $T$ of $G$ such that $S\leq T\leq R$. 
\end{enumerate}
\end{lemma}

\begin{proofof}
To facilitate notation somewhat, set $V \coloneqq [R/S](G/K)$.

\ref{Lem:Seconda}  Let $W=[R/S](G/L)$. Then 
\[
V\big/(L/K) = \left[ (R/S)(L/K)\big/(L/K) \right](G/K)\big/(L/K) \in \F,
\]
where the pairs $\left((R/S)(L/K)\big/(L/K), (G/K)\big/(L/K)\right)$ and $(R/S, G/L)$ are  
equivalent, so $W\cong V\big/(L/K)\in \F$ by Lemma~\ref{Lem:First}.

\ref{Lem:Secondb} Let 
\[
V_0 \coloneqq \left[ (E\cap R)S \big/ S \right](KE/K),
\] 
where $(eS)^{lK} = l^{-1}elS$ for all $eS\in (E\cap R)S/S$ and $lK\in KE/K$ and let
\[
W \coloneqq \left[ (E\cap R) \big/ (E\cap S) \right]\left( E\big/(E\cap K) \right).
\]
Then $V_0 \in \F$ since the class $\F$ is hereditary. 
On the other hand, the pairs 
\[
\left( (E\cap R)S \big/S, KE/K\right), \left( (E\cap R) \big/ (E\cap S), E/(E\cap K) \right)
\]
are equivalent, so $W \cong V_0 \in \F$.

\ref{Lem:Secondc} Observe that $[T/S](G/S)$ is a normal subgroup of
 $[R/S](G/S)$,
thus $[T/S](G/K) \in \F$ since $\F$ is normally hereditary. 
Now let $D \coloneqq [R/T] (G/K)$, where $(rT)^{gK}=r^{g}T$ for all $rT\in R/T$ and 
$gK\in G/K$.
Then  
\[
V \big/ (R/T) = \left[ (R/S) \big/ (T/S) \right]\left( (T/S)(G/K) \big/(T/S) \right) \in \F,
\] 
where the pairs 
\[
(R/T, G/K), \left( (R/S) \big/(T/S), (T/S)(G/K) \big/(T/S) \right)
\] 
are equivalent, so $D \in \F$ by Lemma~\ref{Lem:First}.
The proof is now complete.
\end{proofof}

\begin{lemma}\label{Lem:Third}
Let $N, M$ and $K < H \leq G$ be normal subgroups of $G$.
\begin{enumerate}[label={\upshape(\roman*)}]
\item\label{Lem:Thirda} If $N\leq K$ then
\[
[H/K]\left( G \big/ \cent_{G}(H/K) \right) \cong \left[ (H/N) \big/ (K/N) \right]\left( (G/N)\big/\cent_{G/N}\left( (H/N)\big/(K/N)\right) \right).
\]
\item\label{Lem:Thirdb} If $T/L$ is a normal section of $G$ and $H/K$ and $T/L$ are $G$-isomorphic
then $\cent_{G}(H/K)=\cent_{G}(T/L)$ and $$[H/K](G/\cent_{G}(H/K))\cong [T/L](G/\cent_{G}(T/L)).$$
\item\label{Lem:Thirdc} $[MN/N] (G/\cent_{G}(MN/N))\cong [M/(M\cap N)](G/\cent_{G}(M/(M\cap N)).$
\end{enumerate}
\end{lemma}

\begin{proofof}
\ref{Lem:Thirda} In view of the $G$-isomorphisms $H/K\cong (H/N)\big/(K/N)$ and 
\[
G/\cent_{G}(H/K)\cong (G/N)/(\cent_{G}(H/K)/N),
\] 
the pairs 
\[
\left( H/K, G\big/\cent_{G}(H/K) \right), \left( (H/N)\big/(K/N), (G/N)\big/\cent_{G/N}\left( (H/N)\big(K/N) \right) \right)
\]
are equivalent. Hence Statement~\ref{Lem:Thirda} is a corollary of Lemma~\ref{Lem:First}.

\ref{Lem:Thirdb}  A direct check shows that $C=\cent_{G/N}(H/K) = \cent_{G}(T/L)$ and that the pairs 
$(H/K, G/C)$ and $(T/L, G/C)$ are equivalent. 
Hence Statement~\ref{Lem:Thirdb} is also a corollary of Lemma~\ref{Lem:First}.

\ref{Lem:Thirdc} This follows from the $G$-isomorphism $MN/N\cong M\big/(M\cap N)$ and Statement~\ref{Lem:Thirdb}.
Thus the lemma is proved.
\end{proofof}

We say that a normal subgroup $N$ of $G$ is \emph{$\F$-hypercentral} in $G$ 
if either $N=1$ or every chief factor of $G$ below $N$ is $\F$-central in $G$.

\begin{lemma}\label{Lem:Fourth}
Let $\F$ be a hereditary formation and put $Z \coloneqq \Zent_{\F}(G)$.
Let $A$, $B$ and $N$ be subgroups of $G$, where $N$ is normal in $G$.
Then the following hold.
\begin{enumerate}[label={\upshape(\roman*)}]
\item\label{Lem:Fourtha} $Z$ is $\F$-hypercentral in $G$.
\item\label{Lem:Fourthb} If $N \leq Z$ then $Z/N = \Zent_{\F}(G/N)$.
\item\label{Lem:Fourthc} $\Zent_{\F}(B) \cap A \leq \Zent_{\F}(B \cap A)$.
\item\label{Lem:Fourthd} If $B$ is a normal $\F$-hypercentral subgroup of $G$ then 
$BN \big/ N$ is a normal $\F$-hypercentral subgroup of $G/N$.
\end{enumerate}
\end{lemma}

\begin{proofof}
\ref{Lem:Fourtha} It suffices to consider the case $Z = A_1A_2$, 
where $A_1$ and $A_2$ are normal $\F$-hypercentral subgroups of $G$. 
Moreover, in view of the Jordan-H\"{o}lder theorem for chief series,
it suffices to show that if $A_1 \leq K < H \leq A_1A_2$ then $H/K$ is $\F$-central in $G$. 
But in this case, we have $H = A_1(H \cap A_2)$ and 
\[
H/K = A_1(H \cap A_2) \big/ K \cong_{G} (H\cap A_2) \big/ (K\cap A_2),
\] 
where $(H\cap A_2)\big/ (K\cap A_2)$ is $\F$-central in $G$ and hence $H/K$ is $\F$-central in $G$ by Lemma~\ref{Lem:Third}~\ref{Lem:Thirdb}.

\ref{Lem:Fourthb} This assertion is a corollary of Lemma~\ref{Lem:Third}~\ref{Lem:Thirda}, Statement~\ref{Lem:Fourtha} and the Jordan-H\"{o}lder theorem for chief series.

\ref{Lem:Fourthc} First, assume that $B = G$ and let $1 = Z_0 < Z_1 < \ldots < Z_t \coloneqq Z$ be a chief series of $G$ below $Z$.
Put $C_i \coloneqq \cent_{G}(Z_i / Z_{i-1})$. 
Now  consider the series
\[
1= Z_0 \cap A \leq Z_1 \cap A \leq \ldots  \leq Z_t \cap A = Z \cap A.
\]

Let $i \in \{1,  \ldots , t \}$. Then, by Statement~\ref{Lem:Fourtha},
$Z_{i}/Z_{i-1} $ is $\F$-central in $G$, so $(Z_i \cap A) \big/ (Z_{i-1} \cap A)$
is $\F$-central in $A$ by Lemma~\ref{Lem:Second}~\ref{Lem:Secondb}.
Hence, in view of the Jordan-H\"{o}lder theorem for chief series and Lemma~\ref{Lem:Second}~\ref{Lem:Secondc}, 
we have  $Z \cap A \leq \Zent_{\F}(A)$.

Finally, assume that $B$ is any subgroup of $G$. Then, in view of the preceding paragraph, 
we have
\[
\Zent_{\F}(B) \cap A = \Zent_{\F}(B) \cap (B\cap A)\leq \Zent_{\F}(B \cap A).
\]

\ref{Lem:Fourthd} If $H/K$ is a chief factor of $G$ such that $N \leq K < H \leq NB$, 
then from the $G$-isomorphism 
\[
H/K \cong_{G}(H\cap N) \big/ (K\cap N)
\]
we get that $H/K$ is $\F$-central in $G$ by Lemma~\ref{Lem:Third}~\ref{Lem:Thirdb}, 
so every chief factor of $G/N$ below $BN/N$ is $\F$-central in $G/N$ by Lemma~\ref{Lem:Third}~\ref{Lem:Thirda}. 
Therefore, $BN/N$ is $\F$-hypercentral in $G/N$.  
This concludes the proof of the lemma.
\end{proofof}

\begin{lemma}\label{Lem:Fifth}
Let $\F$ be a saturated formation. 
Then the following statements are equivalent.
\begin{enumerate}[label={\upshape(\roman*)}]
\item\label{Lem:Fiftha} $G\in \F$.
\item\label{Lem:Fifthb} Every chief factor of $G$ is  $\F$-central in $G$.
\item\label{Lem:Fifthc} $G$ has a normal $\F$-hypercentral subgroup $N$ such that $G/N\in \F$.
\end{enumerate}
\end{lemma}
\begin{proofof}
\ref{Lem:Fiftha} $\Rightarrow$ \ref{Lem:Fifthb} This follows directly from the Barnes-Kegel result~\cite{barnes}. 

\ref{Lem:Fifthb} $\Rightarrow$ \ref{Lem:Fifthc} This implication is evident.

\ref{Lem:Fifthc} $\Rightarrow$ \ref{Lem:Fiftha} We prove this implication by induction on $|G|$. 
Let $R$ be a minimal normal subgroup of $G$ and write $C = \cent_{G}(R)$. 
In view of Lemma~\ref{Lem:Fourth}~\ref{Lem:Fourthd}, the hypothesis holds for $G/R$, so $G/R \in \F$ by induction.
Therefore, $G\in \F$ if either $G$ has a minimal normal subgroup $N \ne R$ or if $R \leq \fratt(G)$. 

Now suppose that $R \nleqslant \fratt(G)$ is the unique minimal normal subgroup of $G$ and let 
$M$ be a maximal subgroup of $G$ such that $G = RM$. 
If $R$ is non-abelian, then $C = 1$ since $C$ is normal in $G$ and in this case we have $R \nleqslant C$. 
Therefore, $G \cong G/C\in \F$ by hypothesis. 
Finally, suppose that $R$ is abelian. 
Then $R = C$ since in this case we have $C = R(C\cap M)$, where $C\cap M$ is normal in $G$. 
Therefore, $G = R\rtimes M \in \F$ by Lemma~\ref{Lem:First} since $R/1$ is $\F$-central in $G$ 
and the pairs $(R, M)$ and $(R, G/\cent_R(R)) = (R, G/R)$ are equivalent. 
Our proof is complete.
\end{proofof}

\begin{lemma}\label{Lem:Sixth}
Let $\F$ be a saturated formation and assume that $N$ is normal in $G$.
\begin{enumerate}[label={\upshape(\roman*)}]
\item\label{Lem:Sixtha} If $G/N \in \F$ and $U$ is a minimal supplement to $N$ in $G$ then $U\in \F$.
\item\label{Lem:Sixthb} If $U$ is a subgroup of $G$ such that $U \in \F$ and $NU = G$ 
then $Z \coloneqq U\cap \cent_{G}(N)$ is a normal subgroup of $G$ such that $Z \leq \Zent_{\F}(G)$.
\end{enumerate}
\end{lemma}

\begin{proofof} 
\ref{Lem:Sixtha} This follows from the fact that $U \cap N \leq \fratt(U)$ owing to the minimality of $U$. 

\ref{Lem:Sixthb} Since $G=NU$, $Z$ is normal in $G$. 
Moreover, if $H/K$ is a chief factor of $U$ below $Z$ then $H/K$ is $\F$-central in $U$ by Lemma~\ref{Lem:Fifth} 
and $H/K$ is a chief factor of $G$ since $N\leq \cent_{G}(Z)$. 
Also, we have $N \leq \cent_{G}(H/K)$ and so  
\[
\cent_{G}(H/K) = N \left( \cent_{G}(H/K) \cap U \right) = N\cent_{U}(H/K),
\] 
which in turn implies that  
\begin{align}
G/\cent_{G}(H/K) = NU \big/ N\cent_{U}(H/K) 	&\cong 	U \big/ \left(U\cap N\cent_{U}(H/K) \right)\\
										&=		U\big/ \left(\cent_{U}(H/K)(U\cap N) \right)\\
										&=		U\big/ \cent_{U}(H/K).
\end{align}
Thus the pairs $\left( H/K, U \big/\cent_{U}(H/K) \right)$ and $\left( H/K, G \big/\cent_{G}(H/K) \right)$ are 
equivalent and it follows that $H/K$ is $\F$-central in $G$ by Lemma~\ref{Lem:First}. 
Therefore, $Z \leq \Zent_{\F}(G)$ as wanted.
\end{proofof}

\begin{proofofthmB} 
Let $D = G^{\F}$, $C = \cent_{G}(D)$ and $Z = \Zent_{\F}(G)$. 
Let $U$ be a minimal supplement to $D$ in $G$ and write $T = CU$.

It follows that $T/C \cong U \big/ (U\cap C) \in \F$ by Lemma~\ref{Lem:Sixth}~\ref{Lem:Sixtha}, so $T^{\F} \leq C$. 
On the other hand, we have $T^{\F} \leq D$ since 
\[
T \big/ (T\cap D) \cong TD/D = CUD/D = G/D \in \F
.\]  
Therefore, $T \big/ (C\cap D)\in \F$ and so, in view of Lemma~\ref{Lem:Sixth}~\ref{Lem:Sixtha} again, 
for some subgroup $H$ of $T$ we have $H \in \F$ and $T = (C\cap D)H$. 
Since $C$ is a subgroup of $T$, we have $C = (C \cap D)(C \cap H)$. 
Thus 
\[
G = DU \leq DT = D(C\cap D)H = DH.
\] 
It follows that  $C \cap H \leq Z$ by Lemma~\ref{Lem:Sixth}~\ref{Lem:Sixthb}, 
so $C \leq DZ$ and the theorem is proved.
\end{proofofthmB}

We have now assembled all the necessary tools to begin the proof of our main result.

\begin{proofofthmA}
Assume that the assertion is false and let $G$ be a counterexample with $|G| + |S|$ minimal. 
Then $S < G$ by Theorem~\ref{Thm:thmB}. 
By hypothesis, there is a subgroup chain 
\[
S = S_0 \leq S_1 \leq \ldots \leq S_n = G
\]
such that either $S_{i-1} \trianglelefteq S_i$ or $S_i \big/ (S_{i-1})_{S_i} \in \F$ for all $i \in \{1,\ldots,n\}$.
Since $S$ is a proper subgroup of $G$, we can assume without loss of generality that
 $M \coloneqq S_{n-1} < G$. 
Now write $D = S^{\F}$ and $C = \cent_{G}(D)$, so that $C \nleqslant D$.

\textbf{(Step 1).} Since $\F$ is hereditary and 
\[
S \big/ (S\cap G^{\F}) \cong SG^{\F} \big/G^{\F} \leq G \big/G^{\F},
\] 
we have $D \leq G^{\F}$. Also, from Theorem~\ref{Thm:thmB}  and the hypothesis, we have 
\[
C \cap S = \cent_{S}(D) \leq \Zent_{\F}(S)D = D.
\] 

\textbf{(Step 2).} The hypothesis holds for $(W, S)$ for every proper subgroup $W$ of $G$
 containing $S$ 
by \cite[Lemma 6.1.7]{classes}, so we have $\cent_{W}(D) \leq D$ owing to the choice of $G$. 
Now note that $S \leq \norm_{G}(C)$ since $S \leq \norm_{G}(D)$, 
so $SC$ is a subgroup of $G$ and hence $SC = G$ (otherwise, $C = \cent_{SC}(D) \leq D$). 

\textbf{(Step 3).} Observe that $D \trianglelefteq G$ by virtue of $D$ being characteristic in $S$ and (Step 2).

\textbf{(Step 4).} Let $S \leq M \leq V$, where $V$ is a maximal subgroup of $G$. 
Then $\cent_V(D) \leq D$ by (Step 2), so 
\[
V = V \cap SC = S (V \cap C) = S \cent_V(D) = SD \leq S
\]
and hence $S = V = M$. 
Finally, from $C \cap S \leq D$ it follows that 
\[
G/D = [CD/D](S/D).
\] 

\textbf{(Step 5).} First assume that $S = M$ is not normal in $G$. 
Then $G/S_G \in \F$ and so $D \leq G^{\F} \leq S_G$, where $D$ is normal in $G$ by (Step 2). 
Now (Step 3) implies that 
\[
G/D = [CD/D](S/D),
\]
where $S/D$ is a maximal subgroup of $G/D$ and thus $CD/D$ is a minimal normal subgroup of $G/D$.
Similarly, $CS_G/S_G$ is a minimal normal subgroup of $G/S_G$ owing to 
\[
CS_G\cap S = S_G(C \cap S) \leq S_GD \leq S_G.
\]
Therefore, $CS_G/S_G$ is $\F$-central in $G/S_G$ by Lemma~\ref{Lem:Fifth} and the fact that $G/S_G \in \F$.
 
Then, in view of the $G$-isomorphisms 
\[
CD \big/ C \cong C \big/ (C\cap D) = C \big/ (C\cap S_G) \cong CS_G \big/ S_G,
\] 
the factor $CD/C$ of $G$ is $\F$-central in $G$ by Lemma~\ref{Lem:Third}~\ref{Lem:Thirdb}. 
Therefore, $G/D\in \F$ by Lemma~\ref{Lem:Fifth}, which in turn implies that $D \leq G^{\F} \leq D$ and thus $D = G^{\F}$.

Finally, suppose that  $S$ is normal in $G$ so that $G/S$ is a cyclic group of prime order by (Step 3).  
It follows that 
\[
G/D = [CD/D](S/D) = (CD/D) \times (S/D)
\]
and so  $G/D \in \F$, since $[CD/C](G/C_{G}(CD/D))\in \F$   and $\F$ is hereditary.

\textbf{(Step 6).} From (Step 5) we have $D = G^{\F}$, hence 
\[
C \leq \Zent_{\F}(G)G^{\F} = G^{\F} = D
\] 
by Theorem~\ref{Thm:thmB}, against our assumption on $(G,S)$. 
With this final contradiction our proof is concluded.
\end{proofofthmA}

\begin{proofofcor}
Suppose first that $X$ is a group with no non-trivial normal $\F$-central subgroups. 
In other words, assume that $\Zent_{\F}(X) = 1$. 
In this case, it is a consequence of Theorem~\ref{Thm:thmB}  that 
$\cent_{X}(X^{\F}) \leq X^{\F}$. 
Since $X^{\F}$ is a normal subgroup of $X$, however, it follows that 
$X \big/ X^{\F}$ embeds isomorphically as a subgroup of $\aut(X^{\F})$. 
Thus we see that 
\[
\left\lvert X \big/ X^{\F} \right\rvert \leq \left\lvert X\big/ \cent_{G}(X^{\F}) \right\rvert \leq \left\lvert \aut(X^{\F}) \right\rvert,
\]
and we deduce that $|X| \leq \left\lvert \hol(X^{\F}) \right\rvert$. So the 
assertion holds for groups with trivial $\F$-hypercentre. 

Now, observe that the group $X = G\big/\Zent_{\F}(G)$ has trivial $\F$-hypercentre
 by Lemma~\ref{Lem:Fourth}~\ref{Lem:Fourthb}. 
Also, 
\[
X^{\F} = U\Zent_{\F}(G) \big/ \Zent_{\F}(G) \cong U \big/ (U\cap \Zent_{\F}(G)) \cong U.
\]
The first equality is a consequence of that fact that residuals behave well with respect
 quotients and 
the last isomorphism follows from $U\cap \Zent_{\F}(G) = 1$, which holds by hypothesis. 
Thus
\[
\left\lvert G \big/ \Zent_{\F}(G) \right\rvert \leq \left\lvert \hol \left(U\Zent_{\F}(G)
 \big/ \Zent_{\F}(G) \right) \right\rvert  = \left\lvert \hol(U) \right\rvert,
\] 
as wanted. The corollary is proved.               
\end{proofofcor}

\section{Further applications}\label{Sec:Final}
\textbf{1.} Let $\F$ be a hereditary saturated formation. 
A subgroup $A$ of $G$ is said to be \emph{$\F$-subnormal} in $G$ if there is a subgroup chain 
$A = A_0 \leq A_1 \leq \ldots \leq A_n = G$ such that $A_i / (A_{i-1})_{A_i} \in \F$ for all $i \in \{1, \ldots , n\}$.
It is clear that every $\F$-subnormal subgroup is also $K$-$\F$-subnormal in the group. 
Therefore, we get from Theorem~\ref{Thm:thmA} the following corollaries.

\begin{corollary}\label{Cor:FinalSection1}
Let $\F$ be a hereditary saturated formation. 
Let $S$ be an $\F$-subnormal subgroup of $G$ and suppose that 
$\Zent_{\F}(E) = 1$ for every subgroup $E$ of $G$ such that $S \leq E$. 
Then $\cent_{G}(S^{\F}) \leq S^{\F}$.
\end{corollary}

\begin{corollary}\label{Cor:FinalSection2}
Let $S$ be a $K$-$\mathfrak{U}$-subnormal subgroup of $G$ and suppose that 
$\Zent_{\mathfrak{U}}(E) = 1$ for every subgroup $E$ of $G$ such that $S \leq E$. 
Then $\cent_{G}(S^{\mathfrak{U}}) \leq S^{\mathfrak{U}}$.
\end{corollary}

In this corollary  $S^{\mathfrak{U}}$ is the supersoluble residual of $S$ and 
$\Zent_{\mathfrak{U}}(E)$ is the supersoluble hypercentre of $E$, that is, the product of all 
normal subgroups $N$ of $E$ such that either $N=1$ or $N > 1$ and every 
chief factor of $E$ below $N$ is cyclic.

\textbf{2.} In recent years, the investigations of many authors have focussed on the so-called
$\sigma$-subnormal subgroups 
(see, for example, \cite{s`kibasigma}, \cite{skibatau}, \cite{skibapstgroups}, \cite{skibasublattices19}, \cite{skibasublattices20}, \cite{guoskiba}, \cite{racsam}, \cite{kamornikovfactorised}, \cite{kamornikovyi}, \cite{kamornikovsigma1}). 
We explain this notion below. 

Let $\sigma$ be a partition of the set of all primes $\mathbb{P}$, that is, 
$\sigma =\{\sigma_{i} \mid i\in I \}$, where $\mathbb{P} = \bigcup_{i\in I} \sigma_{i}$ 
and $\sigma_{i} \cap \sigma_{j} = \emptyset$ for all $i \ne j$. 
The group $G$ is said to be \emph{$\sigma$-primary} 
if $G$ is a $\sigma_{i}$-subgroup for some $i$ and \emph{$\sigma$-nilpotent} 
if $G = G_1 \times \cdots \times G_t$ for some $\sigma$-primary groups 
$G_1, \ldots , G_t$. The symbol $\N _{\sigma}$ denotes the class of all
 $\sigma$-nilpotent groups.

Furthermore, a subgroup $A$ of $G$ is said to be \emph{$\sigma$-subnormal} 
in $G$ if there is a subgroup chain 
\[
A = A_0 \leq A_1 \leq \ldots \leq A_n = G
\]
such that either $A_{i-1} \trianglelefteq A_{i}$ or $A_i / (A_{i-1})_{A_i}$ is ${\sigma}$-primary for all $i \in \{1,\ldots,n\}$. 
It is not difficult to show that $A$ is ${\sigma}$-subnormal in $G$ if and only if 
it is $K$-$\N_{\sigma}$-subnormal in $G$. 
The results below are thus also corollaries of Theorem~\ref{Thm:thmA}. 
 
\begin{corollary}[{\cite{dergacheva}}]\label{Cor:FinalSection3}
Let $S$ be a $\sigma$-subnormal subgroup of $G$ and suppose that 
$\Zent_{\sigma}(E) = 1$ for every subgroup $E$ of $G$ such that $S \leq E$. 
Then $\cent_{G}(S^{{\N}_{\sigma}}) \leq S^{{\N}_{\sigma}}$.                           
\end{corollary}

\begin{corollary}\label{Cor:FinalSection4}
Let $S$ be an ${\N}_{\sigma}$-subnormal subgroup of $G$ and suppose that 
$\Zent_{\sigma}(E)=1$ for every subgroup $E$ of $G$ such that $S\leq E$. 
Then $\cent_{G}(S^{{\N}_{\sigma}}) \leq S^{{\N}_{\sigma}}$.
\end{corollary}                         

By $\Zent_{\sigma}(E)$ we denote the $\sigma$-nilpotent hypercentre of $E$ 
which is defined as the product of all normal subgroups $N$ of $E$ such that either 
$N=1$ or $N > 1$ and every chief factor $H/K$ of $E$ below $N$ is \emph{$\sigma$-central} in $E$,
that is, $[H/K](E/\cent_{E}(H/K))$ is $\sigma$-primary.


\end{document}